\documentclass{emsprocart}

\usepackage[T1]{fontenc}
\usepackage[utf8]{inputenc}

\usepackage{tikz}
\usetikzlibrary{decorations.markings}

\usepackage{subfig}

\contact[piotr.sniady@math.uni.wroc.pl]{Piotr Śniady, Institute of
Mathematics, Polish Academy of Sciences, Śniadeckich 8,
00-956 Warszawa, Poland \\
Institute of Mathematics, University of Wrocław, pl.~Grunwaldzki 2/4, \mbox{50-384 Wrocław}, Poland}

\theoremstyle{definition}

\newcommand{\Sym}[1]{\mathfrak{S}(#1)}
\newcommand{\Alt}[1]{\mathfrak{A}(#1)}
\newcommand{\R}{\mathbb{R}}
\newcommand{\Z}{\mathbb{Z}}
\newcommand{\E}{\mathbb{E}}

\DeclareMathOperator{\Tr}{Tr}
\DeclareMathOperator{\Ch}{Ch}

\title[Combinatorics of asymptotic representation theory]{Combinatorics of asymptotic representation theory}

\author[Piotr Śniady]{Piotr Śniady}

\begin{document}

\begin{abstract}
The representation theory of the symmetric groups $\Sym{n}$ is intimately related to combinatorics: combinatorial objects such as Young tableaux and combinatorial algorithms such as Murnaghan-Nakayama rule. In the limit as $n$ tends to infinity, the structure of these combinatorial objects and algorithms becomes complicated and it is hard to extract from them some meaningful answers to asymptotic questions. 
In order to overcome these difficulties, a kind of dual combinatorics of the representation theory of the symmetric groups was initiated in 1990s. We will concentrate on one of its highlights: Kerov polynomials which express characters in terms of, so called, free cumulants. 
\end{abstract}

\begin{classification}
Primary 20C30; 
Secondary 
05E10, 
46L54     
\end{classification}

\begin{keywords}
representations of symmetric groups, Young diagrams, asymptotic representation theory, free cumulants,
Kerov polynomials
\end{keywords}

\maketitle


\vspace{4ex}

\begin{flushright}
\emph{Dedicated to Augustyn Kałuża, my Teacher.}
\end{flushright}

\vspace{1ex}

This note is a guided tour through some selected topics  of \emph{the asymptotic representation theory of the symmetric groups} and its combinatorics. Our guide will be the formula
\begin{equation}
\tag{$\ast$}
\label{eq:main}
 \overbrace{\Ch_5}^{\text{character}} = \overbrace{R_6 + 15R_4 + 5R_2^2 + 8R_2}^{\text{shape}}.
\end{equation}
The forthcoming sections are devoted to an explanation of the cryptic quantities involved here as well as to exploration of the interesting features of this equality.


\section{Representations and characters}

The left-hand side of \eqref{eq:main} is a \emph{character}, a fundamental object in the \emph{representation theory}.
In this section we will briefly review this theory.

\subsection{Example: representation of the symmetric group $\Sym{3}$}
\label{subsec:example1}
Roughly speaking, the subject of the representation theory is the investigation of the ways in which a given abstract group can be realized concretely as a group of matrices. Before we give a formal definition let us have a look on a simple example.

\begin{figure}[t]
\centering
\subfloat[][]{\begin{tikzpicture}
\draw[->,gray] (-2,0) -- (2,0);
\draw[->,gray] (0,-1) -- (0,2);
\draw[thick] (90:1) node[anchor=south west]{1} -- (210:1) node[anchor=north east]{2} -- (330:1) node[anchor=north west]{3} -- cycle ;  
\fill  (90:1) circle(2pt);
\fill (210:1) circle(2pt);
\fill (330:1) circle(2pt);
\label{fig:triangle}
\end{tikzpicture}}
\qquad \qquad
\subfloat[][]{
\begin{tikzpicture}[line join=round,scale=0.5]
\tikzstyle{mypolygonstyle}  = []
\tikzstyle{cube5back} = [very thick,black!20]
\tikzstyle{cube5front} = [very thick,black!black]
\tikzstyle{cube4back}  = [very thick,black!20]
\tikzstyle{cube4front} = [very thick,black]
\tikzstyle{cube3back}  = [very thick,black!20]
\tikzstyle{cube3front} = [very thick,black]
\tikzstyle{cube2back}  = [very thick,black!20]
\tikzstyle{cube2front} = [very thick,black]
\tikzstyle{cube1back}  = [thick,black!10,dashed]
\tikzstyle{cube1front} = [thick,black!50,dashed]
\tikzstyle{dodecahedronback}  = [very thick,black!20]
\tikzstyle{dodecahedronfront} = [ultra thick]
\tikzstyle{mypolygonstyle5} = []
\tikzstyle{mypolygonstyle4} = []
\tikzstyle{mypolygonstyle3} = []
\tikzstyle{mypolygonstyle2} = []
\tikzstyle{mypolygonstyle1} = []
\draw[dodecahedronback,mypolygonstyle](2.893,.063)--(.233,-1.38)--(-2.139,.649)--(-.944,3.346)--(2.166,2.984)--cycle;
\draw[dodecahedronback,mypolygonstyle](-.233,-3.798)--(-2.893,-3.263)--(-4.071,-.514)--(-2.139,.649)--(.233,-1.38)--cycle;
\draw[dodecahedronback,mypolygonstyle](2.139,-3.849)--(-.233,-3.798)--(.233,-1.38)--(2.893,.063)--(4.071,-1.464)--cycle;
\draw[cube1back ,mypolygonstyle1](.233,-1.38)--(4.071,-1.464)--(.944,-3.346)--(-2.893,-3.263)--cycle;
\draw[dodecahedronback,mypolygonstyle](4.071,.514)--(4.071,-1.464)--(2.893,.063)--(2.166,2.984)--(2.893,3.263)--cycle;
\draw[dodecahedronback,mypolygonstyle](2.893,3.263)--(2.166,2.984)--(-.944,3.346)--(-2.139,3.849)--(.233,3.798)--cycle;
\draw[cube1back ,mypolygonstyle1](.233,-1.38)--(-.944,3.346)--(2.893,3.263)--(4.071,-1.464)--cycle;
\draw[dodecahedronback,mypolygonstyle](-2.139,.649)--(-4.071,-.514)--(-4.071,1.464)--(-2.139,3.849)--(-.944,3.346)--cycle;
\draw[cube1back ,mypolygonstyle1](.233,-1.38)--(-2.893,-3.263)--(-4.071,1.464)--(-.944,3.346)--cycle;
\draw[cube1front,mypolygonstyle1](-4.071,1.464)--(-.233,1.38)--(2.893,3.263)--(-.944,3.346)--cycle;
\draw[cube1front,mypolygonstyle1](-4.071,1.464)--(-2.893,-3.263)--(.944,-3.346)--(-.233,1.38)--cycle;
\draw[dodecahedronfront,mypolygonstyle](-.233,-3.798)--(2.139,-3.849)--(.944,-3.346)--(-2.166,-2.984)--(-2.893,-3.263)--cycle;
\draw[dodecahedronfront,mypolygonstyle](-2.893,-3.263)--(-2.166,-2.984)--(-2.893,-.063)--(-4.071,1.464)--(-4.071,-.514)--cycle;
\draw[cube1front,mypolygonstyle1](4.071,-1.464)--(2.893,3.263)--(-.233,1.38)--(.944,-3.346)--cycle;
\draw[dodecahedronfront,mypolygonstyle](.944,-3.346)--(2.139,-3.849)--(4.071,-1.464)--(4.071,.514)--(2.139,-.649)--cycle;
\draw[dodecahedronfront,mypolygonstyle](-4.071,1.464)--(-2.893,-.063)--(-.233,1.38)--(.233,3.798)--(-2.139,3.849)--cycle;
\draw[dodecahedronfront,mypolygonstyle](-.233,1.38)--(2.139,-.649)--(4.071,.514)--(2.893,3.263)--(.233,3.798)--cycle;
\draw[dodecahedronfront,mypolygonstyle](-2.166,-2.984)--(.944,-3.346)--(2.139,-.649)--(-.233,1.38)--(-2.893,-.063)--cycle;
\end{tikzpicture}
\label{fig:dodecahedron}
}

\caption{\protect\subref{fig:triangle}
Equilateral triangle on the plane. 
\protect\subref{fig:dodecahedron}
Regular dodecahedron and one of the five cubes (the dashed lines) which can be inscribed into it.}

\end{figure}
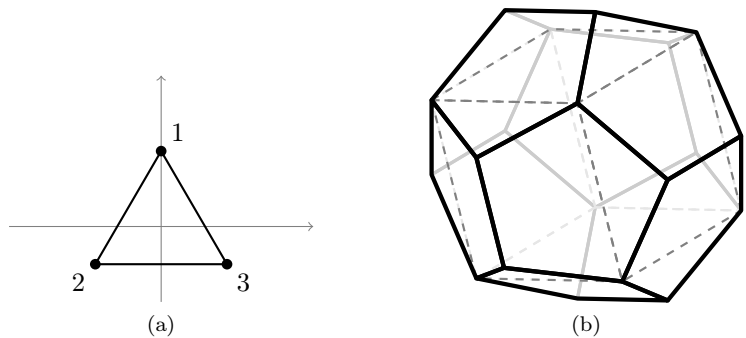

The \emph{symmetric group} $\Sym{3}$ is the group of  permutations of the set $\{1,2,3\}$. If we label the vertices of an equilateral triangle by the elements of this set (Figure \ref{fig:triangle}), any element of $\Sym{3}$ gives rise to a symmetry of the triangle, thus to an isometry of the plane. 
If the coordinate system is chosen properly, these isometries are linear
and thus described by $2\times 2$ matrices. 
We can say (abusing a bit the terminology) that we represented the symmetric group $\Sym{3}$ as certain 
$2\times 2$ matrices.

Formally, a \emph{representation} of a group $G$ is a 
homomorphism $\rho: G \rightarrow M_{n}$ which to the elements of the group associates invertible matrices. 



\subsection{Example: representation of the alternating group $\Alt{5}$}
\label{subsec:dodecahedron}

The above example was too simple; we will now present a less obvious one. 

It is possible to inscribe a cube into the regular dodecahedron in such a way that each vertex of the cube is also a vertex of the dodecahedron (Figure \ref{fig:dodecahedron}). For a fixed dodecahedron there are five such cubes. 
Thus to any rotation of the dodecahedron 
corresponds a permutation of the cubes. This permutation is even 
or, in other words, belongs to the \emph{alternating group} $\Alt{5}$ and this correspondence is bijective.
We can revert the optics: to any element of the alternating group $\Alt{5}$  we associate the corresponding rotation of the dodecahedron. If the coordinate system is chosen properly, this rotation is a linear isometry. In this way we constructed an interesting representation of the alternating group $\Alt{5}$.

\subsection{Motivations}
As we already mentioned, the representation theory studies the ways a given abstract group can be represented concretely as a group of matrices. In this way we can use the power of linear algebra in
order to study problems from the group theory. 

Another motivation comes from harmonic analysis.
One of the most powerful tools for analysis and probability on the real line $\R$ is the Fourier transform. 
If we would like to replace the real line $\R$ by the finite cyclic group $\Z_n$
one should simply use the discrete Fourier transform instead. It is less obvious how to define the
Fourier transform on a \emph{non-commutative} finite group $G$. It turns out that
representations of the group $G$ 
are the right tool to define such an analogue. 


\subsection{Characters}
If we view a representation $\rho:G \to M_n$ as a matrix-valued function, 
its values depend on the choice of the coordinate system in the vector space. Sometimes it would be preferable
to have some quantities which do not depend on such choices. One of such quantities is the \emph{trace} of a matrix. This motivates the study of the \emph{character} of $\rho$
$$ \chi^\rho(g) := \Tr \rho(g) \qquad \text{for } g\in G,$$
which is a scalar-valued function on the group $G$.
A significant part of the representation theory is devoted to investigation of such characters. At first sight it might appear 
that changing the focus from representations to characters might cause some loss of information because
a matrix contains much more data than just its trace. Surprisingly, it is not the case as
almost all natural questions of the representation theory can be reformulated in the language of characters.

The left-hand side of our guiding formula \eqref{eq:main} is such a character (up to some normalizing factors which will be 
discussed later).

\section{Young diagrams and their shapes}

The right-hand side of our favorite equality \eqref{eq:main} describes the \emph{shape} of a \emph{Young diagram}.
This section is devoted to this concept.

\subsection{Irreducible representations}
\label{subsec:irreducible}
If $\rho_1:G\rightarrow M_{n_1}$ and $\rho_2:G\rightarrow M_{n_2}$ are representations of the same group $G$,
we can define a new representation of $G$, called \emph{direct sum} $\rho_1\oplus \rho_2:G\rightarrow M_{n_1+n_2}$ which is given by the block matrices
$$ (\rho_1\oplus\rho_2)(g) = 
\begin{bmatrix}
 \rho_1(g) &    \\
           & \rho_2(g)
\end{bmatrix}. 
$$
Representations which can be written (possibly after change of the coordinate system) as direct sums of smaller representations  are called \emph{reducible} and are less interesting. Our attention will concentrate on representations for which such a decomposition is not possible; they are called \emph{irreducible representations} and  play a fundamental role in the representation theory
(for example they are used in the construction of the non-commutative Fourier transform). The corresponding
characters are called \emph{irreducible characters} and they are in the focus of this article.

\subsection{Irreducible representations of the symmetric groups and Young diagrams}


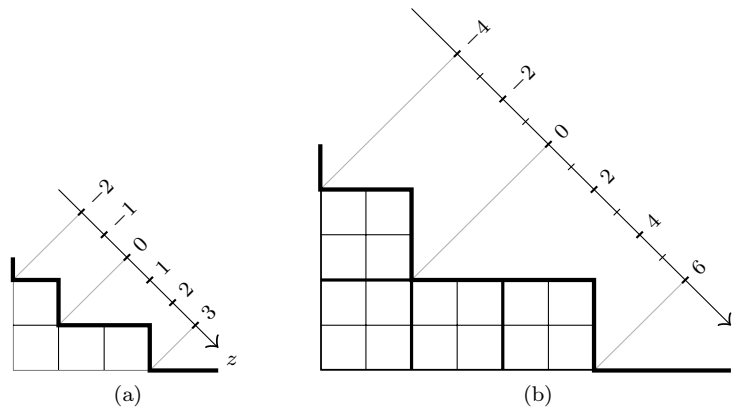
\begin{figure}[t]
\centering
\subfloat[][]{
\begin{tikzpicture}[scale=0.6]
\footnotesize
\draw[black!30](3,0) -- (4,1); 
\draw[black!30](1,1) -- (2.5,2.5); 
\draw[black!30](0,2) -- (1.5,3.5); 

\draw[thick] (2.5,2.5) ++(0,0) +(-2 pt,-2 pt) -- +(2 pt,2 pt) node[anchor=west, rotate=45] {$0$}; 
\draw[thick] (2.5,2.5) ++(0.5,-0.5) +(-2 pt,-2 pt) -- +(2 pt,2 pt) node[anchor=west, rotate=45] {$1$}; 
\draw[thick] (2.5,2.5) ++(1,-1) +(-2 pt,-2 pt) -- +(2 pt,2 pt) node[anchor=west, rotate=45] {$2$}; 
\draw[thick] (2.5,2.5) ++(1.5,-1.5) +(-2 pt,-2 pt) -- +(2 pt,2 pt) node[anchor=west, rotate=45] {$3$}; 
\draw[thick] (2.5,2.5) ++(-0.5,0.5) +(-2 pt,-2 pt) -- +(2 pt,2 pt) node[anchor=west, rotate=45] {$-1$}; 
\draw[thick] (2.5,2.5) ++(-1,1) +(-2 pt,-2 pt) -- +(2 pt,2 pt) node[anchor=west, rotate=45] {$-2$};

\draw[decoration={markings,mark=at position 1 with {\arrow[scale=2]{>}}},
                    postaction={decorate}](1,4) -- (4.5,0.5) node[anchor=north west]{$z$}; 

\draw[ultra thick] (4.5,0) -- (3,0) -- (3,1) -- (1,1) -- (1,2) -- (0,2) -- (0,2.5); 
\clip (0,0) -- (3,0) -- (3,1) -- (1,1) -- (1,2) -- (0,2); 
\draw (0,0) grid (3,3);
\end{tikzpicture}
\label{fig:youngA}
} \hspace{4ex}
\subfloat[][]{
\begin{tikzpicture}[scale=1.2]
\footnotesize
\draw[black!30](3,0) -- (4,1); 
\draw[black!30](1,1) -- (2.5,2.5); 
\draw[black!30](0,2) -- (1.5,3.5); 

\draw[thick] (2.5,2.5) ++(0,0) +(-1 pt,-1 pt) -- +(1 pt,1 pt) node[anchor=west, rotate=45] {$0$}; 
\draw (2.5,2.5) ++(0.25,-0.25) +(-1 pt,-1 pt) -- +(1 pt,1 pt);
\draw[thick] (2.5,2.5) ++(0.5,-0.5) +(-1 pt,-1 pt) -- +(1 pt,1 pt) node[anchor=west, rotate=45] {$2$}; 
\draw (2.5,2.5) ++(0.75,-0.75) +(-1 pt,-1 pt) -- +(1 pt,1 pt);
\draw[thick] (2.5,2.5) ++(1,-1) +(-1 pt,-1 pt) -- +(1 pt,1 pt) node[anchor=west, rotate=45] {$4$}; 
\draw (2.5,2.5) ++(1.25,-1.25) +(-1 pt,-1 pt) -- +(1 pt,1 pt);
\draw[thick] (2.5,2.5) ++(1.5,-1.5) +(-1 pt,-1 pt) -- +(1 pt,1 pt) node[anchor=west, rotate=45] {$6$}; 
\draw (2.5,2.5) ++(-0.25,0.25) +(-1 pt,-1 pt) -- +(1 pt,1 pt);
\draw[thick] (2.5,2.5) ++(-0.5,0.5) +(-1 pt,-1 pt) -- +(1 pt,1 pt) node[anchor=west, rotate=45] {$-2$}; 
\draw (2.5,2.5) ++(-0.75,0.75) +(-1 pt,-1 pt) -- +(1 pt,1 pt);
\draw[thick] (2.5,2.5) ++(-1,1) +(-1 pt,-1 pt) -- +(1 pt,1 pt) node[anchor=west, rotate=45] {$-4$};

\draw[decoration={markings,mark=at position 1 with {\arrow[scale=2]{>}}},
                    postaction={decorate}](1,4) -- (4.5,0.5) node[anchor=north west]{$z$}; 

\draw[ultra thick] (4.5,0) -- (3,0) -- (3,1) -- (1,1) -- (1,2) -- (0,2) -- (0,2.5); 
\clip (0,0) -- (3,0) -- (3,1) -- (1,1) -- (1,2) -- (0,2); 
\draw (0,0) grid[step=0.5] (3,3);
\draw[very thick] (0,0) grid (3,3);
\end{tikzpicture}
\label{fig:youngB}
}

\caption{\protect\subref{fig:youngA} Young diagram $\lambda=(3,1)$ corresponding to the partition $4=3+1$. 
\protect\subref{fig:youngB}
The dilation $2\lambda$ of the Young diagram $\lambda$ shown on the left.}
\end{figure}

There is a bijection between irreducible representations of
the symmetric group $\Sym{n}$ and \emph{Young diagrams} with $n$ boxes. The latter are 
collections of boxes which are nicely aligned to the left and
to the bottom (Figure \ref{fig:youngA}). For a Young diagram $\lambda$ we will denote the corresponding irreducible
representation by $\rho^\lambda$.

Unfortunately, the details of this bijection are technically involved. In order to give 
the flavor of this difficulty we mention only that one of the irreducible representations of the symmetric group
$\Sym{5}$ is closely related to the not-so-trivial representation of its subgroup $\Alt{5}\subset\Sym{5}$ which we discussed in Section \ref{subsec:dodecahedron}.

\subsection{Shape of the Young diagram}

\emph{What can we say about the irreducible representations
of the symmetric groups when the corresponding Young diagrams tend infinity, having a fixed `macroscopic shape'?}
In order to make this question more concrete, we will use the notion of \emph{dilation}.
If $s$ is a positive integer and $\lambda$ is a Young diagram we will denote by $s\lambda$ the 
\emph{dilated Young diagram}, obtained from $\lambda$ by replacing each box by a $s\times s$ grid of boxes (Figure \ref{fig:youngB}). 
If we disregard the size, such a dilated Young diagram has the same \emph{shape}
as the original diagram (compare Figures \ref{fig:youngA} and \ref{fig:youngB}). 
Our original question can be therefore reformulated as investigation of the dilated Young diagrams $s\lambda$, where $\lambda$ is fixed and $s\to\infty$.

\subsection{Homogeneous functions} 
\label{subsec:homogeneous}
For this kind of asymptotic problems we need the right tools:
functions on the set of Young diagrams which would depend \emph{`nicely'} on the shape of the Young diagram.  
For example we could require from such a nice function $f$ that it depends only
on the shape of the Young diagram and not on its size: $f(s \lambda) = f(\lambda)$.
This requirement is too strong; it would not be a big problem if we allow a simple dependence of $f$ on the the size of the Young diagram:
$$ f(s\lambda) = s^k f(\lambda) $$
for some exponent $k$. If this is the case we say that $f$ is \emph{homogeneous} of degree $k$. 

\section{Relationship between characters and the shape?}

Our favorite formula \eqref{eq:main} gives a relationship between the irreducible characters of the symmetric groups
and the shape of the Young diagram. In this section we will investigate this kind of relationships.


\begin{figure}[t]
\centering
\begin{tikzpicture}[scale=0.5]
\scriptsize
\draw[](15,0) -- (0,0) -- (0,10); 
\draw[draw=none,fill=black!0](15.9,0) -- (15,0) -- (15,5) -- (5,5) -- (5,10) -- (0,10) -- (0,10.9) -- (0,0) -- cycle ;
\draw(15.9,0) -- (15,0) -- (15,5) -- (5,5) -- (5,10) -- (0,10) -- (0,10.9); 

\draw[draw=none,fill=black!10](14.9,0) -- (14,0) -- (14,1) -- (12,1) -- (12,2) -- (11,2) -- (11,4) -- (6,4) -- (6,5) -- (4,5) -- (4,6) -- (3,6) -- (3,8) -- (2,8) -- (2,9) -- (0,9) -- (0,9.9) -- (0,0) -- cycle ;
\draw[ultra thick](14.9,0) -- (14,0) -- (14,1) -- (12,1) -- (12,2) -- (11,2) -- (11,4) -- (6,4) -- (6,5) -- (4,5) -- (4,6) -- (3,6) -- (3,8) -- (2,8) -- (2,9) -- (0,9) -- (0,9.9); 

\draw[draw=none,fill=black!30](10.9,0) -- (10,0) -- (10,1) -- (7,1) -- 
(7,2) -- (6,2) -- (6,3)  -- (5,3) -- (5,4) -- (2,4) -- (2,6) -- (1,6) -- (1,7) -- (0,7) -- (0,7.9) -- (0,0) -- cycle ;
\draw[very thick](10.9,0) -- (10,0) -- (10,1) -- (7,1) -- (7,2) -- (6,2) -- (6,3)  -- (5,3) -- (5,4) -- (2,4) -- (2,6) -- (1,6) -- (1,7) -- (0,7) -- (0,7.9); 

\begin{scope}
\clip(15,0) -- (15,5) -- (5,5) -- (5,10) -- (0,10) -- (0,0);
\draw (0,0) grid (15,10);
\end{scope}

\draw (0.5,0.5) node {$1$};
\draw (0.5,1.5) node {$2$};
\draw (1.5,0.5) node {$3$};
\draw (2.5,0.5) node {$4$};
\draw (3.5,0.5) node {$5$};
\draw (4.5,0.5) node {$6$};
\draw (1.5,1.5) node {$7$};
\draw (0.5,2.5) node {$8$};
\draw (2.5,1.5) node {$9$};
\draw (1.5,2.5) node {$10$};
\draw (0.5,3.5) node {$11$};
\draw (0.5,4.5) node {$12$};
\draw (2.5,2.5) node {$13$};
\draw (5.5,0.5) node {$14$};
\draw (3.5,1.5) node {$15$};
\draw (4.5,1.5) node {$16$};
\draw (1.5,3.5) node {$17$};
\draw (0.5,5.5) node {$18$};
\draw (2.5,3.5) node {$19$};
\draw (1.5,4.5) node {$20$};
\draw (3.5,2.5) node {$21$};
\draw (3.5,3.5) node {$22$};
\draw (4.5,2.5) node {$23$};
\draw (5.5,1.5) node {$24$};
\draw (6.5,0.5) node {$25$};
\draw (0.5,6.5) node {$26$};
\draw (6.5,1.5) node {$27$};
\draw (7.5,0.5) node {$28$};
\draw (5.5,2.5) node {$29$};
\draw (4.5,3.5) node {$30$};
\draw (8.5,0.5) node {$31$};
\draw (9.5,0.5) node {$32$};
\draw (1.5,5.5) node {$33$};
\draw (6.5,2.5) node {$34$};
\draw (2.5,4.5) node {$35$};
\draw (3.5,4.5) node {$36$};
\draw (2.5,5.5) node {$37$};
\draw (1.5,6.5) node {$38$};
\draw (7.5,1.5) node {$39$};
\draw (10.5,0.5) node {$40$};
\draw (8.5,1.5) node {$41$};
\draw (4.5,4.5) node {$42$};
\draw (5.5,3.5) node {$43$};
\draw (9.5,1.5) node {$44$};
\draw (7.5,2.5) node {$45$};
\draw (5.5,4.5) node {$46$};
\draw (8.5,2.5) node {$47$};
\draw (10.5,1.5) node {$48$};
\draw (9.5,2.5) node {$49$};
\draw (11.5,0.5) node {$50$};
\draw (0.5,7.5) node {$51$};
\draw (6.5,3.5) node {$52$};
\draw (12.5,0.5) node {$53$};
\draw (2.5,6.5) node {$54$};
\draw (7.5,3.5) node {$55$};
\draw (1.5,7.5) node {$56$};
\draw (11.5,1.5) node {$57$};
\draw (0.5,8.5) node {$58$};
\draw (3.5,5.5) node {$59$};
\draw (1.5,8.5) node {$60$};
\draw (13.5,0.5) node {$61$};
\draw (2.5,7.5) node {$62$};
\draw (10.5,2.5) node {$63$};
\draw (8.5,3.5) node {$64$};
\draw (9.5,3.5) node {$65$};
\draw (10.5,3.5) node {$66$};
\draw (6.5,4.5) node {$67$};
\draw (7.5,4.5) node {$68$};
\draw (12.5,1.5) node {$69$};
\draw (8.5,4.5) node {$70$};
\draw (11.5,2.5) node {$71$};
\draw (2.5,8.5) node {$72$};
\draw (14.5,0.5) node {$73$};
\draw (11.5,3.5) node {$74$};
\draw (0.5,9.5) node {$75$};
\draw (12.5,2.5) node {$76$};
\draw (13.5,1.5) node {$77$};
\draw (9.5,4.5) node {$78$};
\draw (3.5,6.5) node {$79$};
\draw (13.5,2.5) node {$80$};
\draw (1.5,9.5) node {$81$};
\draw (10.5,4.5) node {$82$};
\draw (12.5,3.5) node {$83$};
\draw (11.5,4.5) node {$84$};
\draw (13.5,3.5) node {$85$};
\draw (14.5,1.5) node {$86$};
\draw (4.5,5.5) node {$87$};
\draw (12.5,4.5) node {$88$};
\draw (2.5,9.5) node {$89$};
\draw (13.5,4.5) node {$90$};
\draw (14.5,2.5) node {$91$};
\draw (4.5,6.5) node {$92$};
\draw (3.5,7.5) node {$93$};
\draw (3.5,8.5) node {$94$};
\draw (4.5,7.5) node {$95$};
\draw (14.5,3.5) node {$96$};
\draw (14.5,4.5) node {$97$};
\draw (3.5,9.5) node {$98$};
\draw (4.5,8.5) node {$99$};
\draw (4.5,9.5) node {$100$};
\end{tikzpicture}
\caption{Example of a Young tableau. Shaded regions show boxes with numbers smaller than some thresholds.}
\label{fig:tableau}
\end{figure}

For a wide class of questions concerning irreducible representations of the symmetric 
groups there is a known answer given in terms of some combinatorial algorithm involving boxes of the Young diagram.
For example, the dimension of the irreducible representation $\rho^\lambda$ is equal to the number of \emph{Young 
tableaux} filling $\lambda$. A Young tableau is a filling of the boxes of the Young diagram $\lambda$ with numbers 
$1,2,\dots,n$ (where $n$ is the number of boxes of $\lambda$) in such a way that each number is used exactly once and
the numbers increase from left to right and from bottom to top (Figure \ref{fig:tableau}). Investigation of 
algorithms with a similar flavor is a one of important branches of combinatorics.

In particular, irreducible characters of the symmetric groups
$$ \chi^\lambda(\pi) := \Tr \rho^\lambda(\pi) \qquad \text{for } \pi\in \Sym{n}$$
can be calculated using such a combinatorial algorithm, the \emph{Murnaghan-Nakayama rule} which is a signed sum over, roughly speaking, Young tableaux filling $\lambda$ with some additional properties (related to the conjugacy class of the permutation $\pi$).

Unfortunately, it is a common feature of such combinatorial algorithms that they quickly become
cumbersome when the number of the boxes of the Young diagram tends to infinity. For this reason they
are not very suitable for the investigation of the asymptotics problems. In particular, Murnaghan-Nakayama rule
does not give too much insight into our favorite question --- the answer to which is given by our guiding
equality \eqref{eq:main} --- about the relationship between the characters and the shape of the Young diagram.
In order to overcome these difficulties we will have to find a better normalization of the characters
as well as find a good way of describing the shape of a Young diagram.

\section{Normalized characters}

The left-hand side of our favorite equality \eqref{eq:main} is the \emph{normalized character}.
In this section we will present the details of this quantity.

The usual way of studying the characters 
of the symmetric groups is to fix the Young diagram $\lambda$ and to view $\chi^\lambda(\pi)$ as a function
of the permutation $\pi$. It was a brilliant idea of Kerov and Olshanski to do the opposite and to study
the \emph{dual combinatorics} of representations of the symmetric groups, see below.

For a fixed integer $k\geq 1$ we will denote by 
$[k]=(1,2,\dots,k)\in \Sym{k}$ the full cycle; we will investigate the characters evaluated on the permutation
$[k]$. 
Let $\lambda$ be a Young diagram with $n$ boxes; we are interested in the character
$$ \Tr \rho^{\lambda}([k]). $$
It might seem that this quantity does not make much sense since $[k]$ belongs to $\Sym{k}$ while $\rho^\lambda$ is a representation of a different group, namely $\Sym{n}$. Nevertheless, for $n\geq k$ we can consider an embedding
$\Sym{k}\subset\Sym{n}$ and regard $[k]$ as an element of $\Sym{n}$ simply by adding additional fixpoints.

%
%

It turns out that the `right' way to define the \emph{normalized character} on a cycle of length $k$ is as follows:
$$ \Ch_k(\lambda) := \underbrace{n (n-1) \cdots (n-k+1)}_{\text{$k$ factors}} \frac{\Tr \rho^\lambda([k])}{\Tr \rho^\lambda(e)}, $$
where $n$ is the number of boxes of $\lambda$. As we already mentioned, we will view the normalized character
$\Ch_k$ as a function on the set of Young diagrams and \emph{we impose no restrictions on the
number of boxes of the Young diagrams}.

The normalization factor in the above definition 
is equal to zero if $n<k$; in this way the right-hand side is equal to zero and we do not have to worry that $\rho^\lambda([k])$ is not well-defined in this case. 
The denominator $\Tr \rho^\lambda(e)$, the character on the group unit, is equal just to the dimension of the representation $\rho^\lambda$ and there are some effective methods of calculating it.
This means that the normalized characters $\Ch_k$ contain essentially the same information as the usual characters 
$\chi^\lambda(\pi)$ thus they are just as interesting.

On the other hand the normalized characters $\Ch_k$ have some advantages over the usual characters $\chi^\lambda(\pi)$,
for example
$$ (\lambda_1,\lambda_2,\dots) \mapsto \Ch_k(\lambda) $$
is a polynomial function of the lengths of the rows of the Young diagram $\lambda=(\lambda_1,\lambda_2,\dots)$.



\section{Free cumulants}

The right-hand side of our favorite formula \eqref{eq:main} concerns the shape of the Young diagram.
The question is: \emph{how to choose parameters which would describe the shape of the Young diagram
in the best way?}  The answer comes from Voiculescu's \emph{free probability} theory.

\subsection{Random matrices and free cumulants}
\label{subsec:random-matrices}

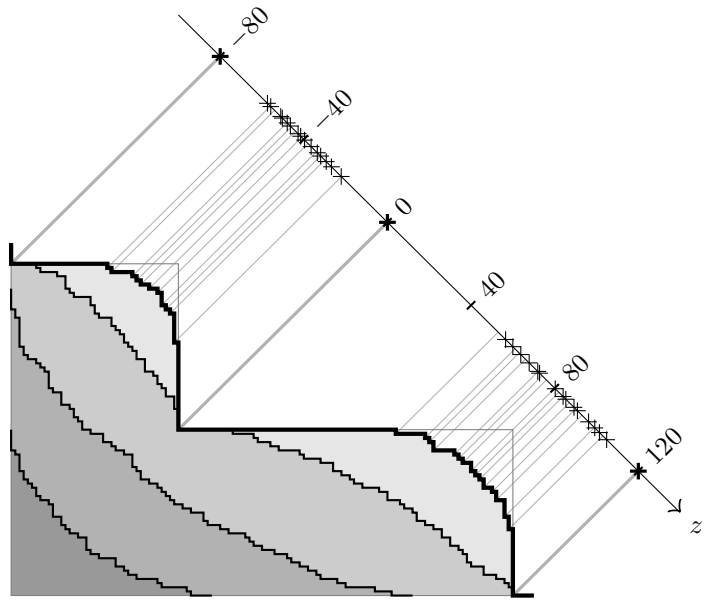
\begin{figure}[t]
\centering
\begin{tikzpicture}[scale=0.055]
\draw[black!30](119,16) -- (141.5,38.5); 
\draw[black!30](118,19) -- (139.5,40.5); 
\draw[black!30](116,23) -- (136.5,43.5); 
\draw[black!30](115,25) -- (135.0,45.0); 
\draw[black!30](113,26) -- (133.5,46.5); 
\draw[black!30](112,28) -- (132.0,48.0); 
\draw[black!30](111,29) -- (131.0,49.0); 
\draw[black!30](110,30) -- (130.0,50.0); 
\draw[black!30](108,32) -- (128.0,52.0); 
\draw[black!30](107,33) -- (127.0,53.0); 
\draw[black!30](106,34) -- (126.0,54.0); 
\draw[black!30](101,35) -- (123.0,57.0); 
\draw[black!30](100,37) -- (121.5,58.5); 
\draw[black!30](99,38) -- (120.5,59.5); 
\draw[black!30](92,39) -- (116.5,63.5); 

\draw[black!30](39,61) -- (79.0,101.0); 
\draw[black!30](38,68) -- (75.0,105.0); 
\draw[black!30](37,69) -- (74.0,106.0); 
\draw[black!30](36,70) -- (73.0,107.0); 
\draw[black!30](35,73) -- (71.0,109.0); 
\draw[black!30](33,74) -- (69.5,110.5); 
\draw[black!30](31,75) -- (68.0,112.0); 
\draw[black!30](30,76) -- (67.0,113.0); 
\draw[black!30](29,77) -- (66.0,114.0); 
\draw[black!30](24,78) -- (63.0,117.0); 
\draw[black!30](23,79) -- (62.0,118.0); 

\draw[black!30,very thick](120,0) -- (150.0,30.0); 
\draw[black!30,very thick](40,40) -- (90.0,90.0); 
\draw[black!30,very thick](0,80) -- (50.0,130.0);

\draw[decoration={markings,mark=at position 1 with {\arrow[scale=2]{>}}},
                    postaction={decorate}](40.0,140.0) -- (160.0,20.0) node[anchor=north west]{$z$}; 

\draw[draw=none,fill=black!10](120,0) -- (120,16) -- (119,16) -- (119,19) -- (118,19) -- (118,23) -- (116,23) -- (116,25) -- (115,25) -- (115,26) -- (113,26) -- (113,28) -- (112,28) -- (112,29) -- (111,29) -- (111,30) -- (110,30) -- (110,32) -- (108,32) -- (108,33) -- (107,33) -- (107,34) -- (106,34) -- (106,35) -- (101,35) -- (101,37) -- (100,37) -- (100,38) -- (99,38) -- (99,39) -- (92,39) -- (92,40) -- (40,40) -- (40,61) -- (39,61) -- (39,68) -- (38,68) -- (38,69) -- (37,69) -- (37,70) -- (36,70) -- (36,73) -- (35,73) -- (35,74) -- (33,74) -- (33,75) -- (31,75) -- (31,76) -- (30,76) -- (30,77) -- (29,77) -- (29,78) -- (24,78) -- (24,79) -- (23,79) -- (23,80) -- (0,80) -- (0,0) -- cycle; 

\draw[draw=none,fill=black!20](120,0) -- (120,2) -- (119,2) -- (119,3) -- (118,3) -- (118,4) -- (116,4) -- (116,5) -- (113,5) -- (113,7) -- (110,7) -- (110,8) -- (109,8) -- (109,10) -- (108,10) -- (108,11) -- (107,11) -- (107,12) -- (105,12) -- (105,14) -- (103,14) -- (103,16) -- (101,16) -- (101,17) -- (99,17) -- (99,18) -- (97,18) -- (97,20) -- (96,20) -- (96,21) -- (94,21) -- (94,22) -- (93,22) -- (93,23) -- (90,23) -- (90,25) -- (88,25) -- (88,26) -- (84,26) -- (84,27) -- (82,27) -- (82,28) -- (81,28) -- (81,29) -- (79,29) -- (79,30) -- (77,30) -- (77,32) -- (74,32) -- (74,33) -- (70,33) -- (70,34) -- (69,34) -- (69,35) -- (64,35) -- (64,36) -- (62,36) -- (62,37) -- (59,37) -- (59,38) -- (58,38) -- (58,39) -- (53,39) -- (53,40) -- (40,40) -- (40,45) -- (39,45) -- (39,48) -- (37,48) -- (37,50) -- (36,50) -- (36,51) -- (35,51) -- (35,53) -- (34,53) -- (34,54) -- (32,54) -- (32,57) -- (30,57) -- (30,59) -- (29,59) -- (29,60) -- (28,60) -- (28,61) -- (27,61) -- (27,63) -- (26,63) -- (26,64) -- (25,64) -- (25,65) -- (24,65) -- (24,67) -- (21,67) -- (21,68) -- (20,68) -- (20,70) -- (19,70) -- (19,72) -- (15,72) -- (15,73) -- (14,73) -- (14,74) -- (13,74) -- (13,77) -- (11,77) -- (11,78) -- (7,78) -- (7,79) -- (6,79) -- (6,80) -- (0,80) -- (0,0) -- cycle ;

\draw[draw=none,fill=black!30](91,0) -- (91,1) -- (87,1) -- (87,2) -- (80,2) -- (80,3) -- (78,3) -- (78,4) -- (74,4) -- (74,5) -- (73,5) -- (73,6) -- (71,6) -- (71,7) -- (69,7) -- (69,8) -- (67,8) -- (67,9) -- (66,9) -- (66,10) -- (63,10) -- (63,11) -- (60,11) -- (60,12) -- (58,12) -- (58,13) -- (57,13) -- (57,14) -- (56,14) -- (56,15) -- (52,15) -- (52,16) -- (50,16) -- (50,17) -- (47,17) -- (47,18) -- (46,18) -- (46,20) -- (43,20) -- (43,21) -- (41,21) -- (41,22) -- (39,22) -- (39,23) -- (38,23) -- (38,24) -- (37,24) -- (37,25) -- (36,25) -- (36,26) -- (35,26) -- (35,27) -- (33,27) -- (33,28) -- (31,28) -- (31,30) -- (30,30) -- (30,31) -- (29,31) -- (29,34) -- (27,34) -- (27,35) -- (25,35) -- (25,36) -- (24,36) -- (24,37) -- (22,37) -- (22,39) -- (19,39) -- (19,40) -- (17,40) -- (17,42) -- (16,42) -- (16,43) -- (15,43) -- (15,45) -- (14,45) -- (14,46) -- (12,46) -- (12,49) -- (10,49) -- (10,50) -- (9,50) -- (9,54) -- (8,54) -- (8,55) -- (7,55) -- (7,56) -- (6,56) -- (6,57) -- (4,57) -- (4,58) -- (3,58) -- (3,60) -- (2,60) -- (2,67) -- (1,67) -- (1,69) -- (0,69) -- (0,0) -- cycle ;

\draw[draw=none,fill=black!40](43,0) -- (43,1) -- (42,1) -- (42,2) -- (36,2) -- (36,3) -- (35,3) -- (35,4) -- (30,4) -- (30,5) -- (28,5) -- (28,6) -- (27,6) -- (27,7) -- (25,7) -- (25,9) -- (21,9) -- (21,10) -- (19,10) -- (19,11) -- (18,11) -- (18,14) -- (16,14) -- (16,15) -- (14,15) -- (14,16) -- (13,16) -- (13,17) -- (12,17) -- (12,19) -- (10,19) -- (10,22) -- (9,22) -- (9,23) -- (7,23) -- (7,24) -- (6,24) -- (6,26) -- (5,26) -- (5,27) -- (4,27) -- (4,28) -- (3,28) -- (3,30) -- (2,30) -- (2,34) -- (1,34) -- (1,35) -- (0,35) -- (0,0) -- cycle ;

\draw[black!50](120,0) -- (120,40) -- (40,40) -- (40,80) -- (0,80) -- (0,0) -- cycle ;

\draw[black,ultra thick](125,0) -- (120,0) -- (120,16) -- (119,16) -- (119,19) -- (118,19) -- (118,23) -- (116,23) -- (116,25) -- (115,25) -- (115,26) -- (113,26) -- (113,28) -- (112,28) -- (112,29) -- (111,29) -- (111,30) -- (110,30) -- (110,32) -- (108,32) -- (108,33) -- (107,33) -- (107,34) -- (106,34) -- (106,35) -- (101,35) -- (101,37) -- (100,37) -- (100,38) -- (99,38) -- (99,39) -- (92,39) -- (92,40) -- (40,40) -- (40,61) -- (39,61) -- (39,68) -- (38,68) -- (38,69) -- (37,69) -- (37,70) -- (36,70) -- (36,73) -- (35,73) -- (35,74) -- (33,74) -- (33,75) -- (31,75) -- (31,76) -- (30,76) -- (30,77) -- (29,77) -- (29,78) -- (24,78) -- (24,79) -- (23,79) -- (23,80) -- (0,80) -- (0,85); 
\draw[thick](125,0) -- (120,0) -- (120,2) -- (119,2) -- (119,3) -- (118,3) -- (118,4) -- (116,4) -- (116,5) -- (113,5) -- (113,7) -- (110,7) -- (110,8) -- (109,8) -- (109,10) -- (108,10) -- (108,11) -- (107,11) -- (107,12) -- (105,12) -- (105,14) -- (103,14) -- (103,16) -- (101,16) -- (101,17) -- (99,17) -- (99,18) -- (97,18) -- (97,20) -- (96,20) -- (96,21) -- (94,21) -- (94,22) -- (93,22) -- (93,23) -- (90,23) -- (90,25) -- (88,25) -- (88,26) -- (84,26) -- (84,27) -- (82,27) -- (82,28) -- (81,28) -- (81,29) -- (79,29) -- (79,30) -- (77,30) -- (77,32) -- (74,32) -- (74,33) -- (70,33) -- (70,34) -- (69,34) -- (69,35) -- (64,35) -- (64,36) -- (62,36) -- (62,37) -- (59,37) -- (59,38) -- (58,38) -- (58,39) -- (53,39) -- (53,40) -- (40,40) -- (40,45) -- (39,45) -- (39,48) -- (37,48) -- (37,50) -- (36,50) -- (36,51) -- (35,51) -- (35,53) -- (34,53) -- (34,54) -- (32,54) -- (32,57) -- (30,57) -- (30,59) -- (29,59) -- (29,60) -- (28,60) -- (28,61) -- (27,61) -- (27,63) -- (26,63) -- (26,64) -- (25,64) -- (25,65) -- (24,65) -- (24,67) -- (21,67) -- (21,68) -- (20,68) -- (20,70) -- (19,70) -- (19,72) -- (15,72) -- (15,73) -- (14,73) -- (14,74) -- (13,74) -- (13,77) -- (11,77) -- (11,78) -- (7,78) -- (7,79) -- (6,79) -- (6,80) -- (0,80) -- (0,85); 
\draw[thick](96,0) -- (91,0) -- (91,1) -- (87,1) -- (87,2) -- (80,2) -- (80,3) -- (78,3) -- (78,4) -- (74,4) -- (74,5) -- (73,5) -- (73,6) -- (71,6) -- (71,7) -- (69,7) -- (69,8) -- (67,8) -- (67,9) -- (66,9) -- (66,10) -- (63,10) -- (63,11) -- (60,11) -- (60,12) -- (58,12) -- (58,13) -- (57,13) -- (57,14) -- (56,14) -- (56,15) -- (52,15) -- (52,16) -- (50,16) -- (50,17) -- (47,17) -- (47,18) -- (46,18) -- (46,20) -- (43,20) -- (43,21) -- (41,21) -- (41,22) -- (39,22) -- (39,23) -- (38,23) -- (38,24) -- (37,24) -- (37,25) -- (36,25) -- (36,26) -- (35,26) -- (35,27) -- (33,27) -- (33,28) -- (31,28) -- (31,30) -- (30,30) -- (30,31) -- (29,31) -- (29,34) -- (27,34) -- (27,35) -- (25,35) -- (25,36) -- (24,36) -- (24,37) -- (22,37) -- (22,39) -- (19,39) -- (19,40) -- (17,40) -- (17,42) -- (16,42) -- (16,43) -- (15,43) -- (15,45) -- (14,45) -- (14,46) -- (12,46) -- (12,49) -- (10,49) -- (10,50) -- (9,50) -- (9,54) -- (8,54) -- (8,55) -- (7,55) -- (7,56) -- (6,56) -- (6,57) -- (4,57) -- (4,58) -- (3,58) -- (3,60) -- (2,60) -- (2,67) -- (1,67) -- (1,69) -- (0,69) -- (0,74); 
\draw[thick](48,0) -- (43,0) -- (43,1) -- (42,1) -- (42,2) -- (36,2) -- (36,3) -- (35,3) -- (35,4) -- (30,4) -- (30,5) -- (28,5) -- (28,6) -- (27,6) -- (27,7) -- (25,7) -- (25,9) -- (21,9) -- (21,10) -- (19,10) -- (19,11) -- (18,11) -- (18,14) -- (16,14) -- (16,15) -- (14,15) -- (14,16) -- (13,16) -- (13,17) -- (12,17) -- (12,19) -- (10,19) -- (10,22) -- (9,22) -- (9,23) -- (7,23) -- (7,24) -- (6,24) -- (6,26) -- (5,26) -- (5,27) -- (4,27) -- (4,28) -- (3,28) -- (3,30) -- (2,30) -- (2,34) -- (1,34) -- (1,35) -- (0,35) -- (0,40); 

\draw (61.304,118.696) +(2,0) -- +(-2,0) +(0,2)--+(0,-2);
\draw (62.092,117.908) +(2,0) -- +(-2,0) +(0,2)--+(0,-2);
\draw (64.486,115.514) +(2,0) -- +(-2,0) +(0,2)--+(0,-2);
\draw (64.858,115.142) +(2,0) -- +(-2,0) +(0,2)--+(0,-2);
\draw (66.032,113.968) +(2,0) -- +(-2,0) +(0,2)--+(0,-2);
\draw (66.786,113.214) +(2,0) -- +(-2,0) +(0,2)--+(0,-2);
\draw (68.632,111.368) +(2,0) -- +(-2,0) +(0,2)--+(0,-2);
\draw (69.242,110.758) +(2,0) -- +(-2,0) +(0,2)--+(0,-2);
\draw (70.1774,109.8226) +(2,0) -- +(-2,0) +(0,2)--+(0,-2);
\draw (71.743,108.257) +(2,0) -- +(-2,0) +(0,2)--+(0,-2);
\draw (73.227,106.773) +(2,0) -- +(-2,0) +(0,2)--+(0,-2);
\draw (74.029,105.971) +(2,0) -- +(-2,0) +(0,2)--+(0,-2);
\draw (75.362,104.638) +(2,0) -- +(-2,0) +(0,2)--+(0,-2);
\draw (76.611,103.389) +(2,0) -- +(-2,0) +(0,2)--+(0,-2);
\draw (78.9482,101.0518) +(2,0) -- +(-2,0) +(0,2)--+(0,-2);
\draw (118.256,61.744) +(2,0) -- +(-2,0) +(0,2)--+(0,-2);
\draw (119.988,60.012) +(2,0) -- +(-2,0) +(0,2)--+(0,-2);
\draw (121.886,58.114) +(2,0) -- +(-2,0) +(0,2)--+(0,-2);
\draw (123.896,56.104) +(2,0) -- +(-2,0) +(0,2)--+(0,-2);
\draw (125.954,54.046) +(2,0) -- +(-2,0) +(0,2)--+(0,-2);
\draw (126.382,53.618) +(2,0) -- +(-2,0) +(0,2)--+(0,-2);
\draw (130.128,49.872) +(2,0) -- +(-2,0) +(0,2)--+(0,-2);
\draw (132.006,47.994) +(2,0) -- +(-2,0) +(0,2)--+(0,-2);
\draw (132.662,47.338) +(2,0) -- +(-2,0) +(0,2)--+(0,-2);
\draw (134.614,45.386) +(2,0) -- +(-2,0) +(0,2)--+(0,-2);
\draw (135.434,44.566) +(2,0) -- +(-2,0) +(0,2)--+(0,-2);
\draw (138.1,41.9) +(2,0) -- +(-2,0) +(0,2)--+(0,-2);
\draw (139.58,40.42) +(2,0) -- +(-2,0) +(0,2)--+(0,-2);
\draw (140.7,39.3) +(2,0) -- +(-2,0) +(0,2)--+(0,-2);
\draw (142.42,37.58) +(2,0) -- +(-2,0) +(0,2)--+(0,-2);

\draw[very thick] (50.0,130.0) +(2,0) -- +(-2,0) +(0,2)--+(0,-2);
\draw[very thick] (90.0,90.0) +(2,0) -- +(-2,0) +(0,2)--+(0,-2);
\draw[very thick] (150.0,30.0) +(2,0) -- +(-2,0) +(0,2)--+(0,-2);

\draw[thick] (90,90) ++(-40,40) +(-1,-1) -- +(1,1) node[anchor=west, rotate=45] {$-80$}; 
\draw[thick] (90,90) ++(-20,20) +(-1,-1) -- +(1,1) node[anchor=west, rotate=45] {$-40$}; 
\draw[thick] (90,90) +(-1,-1) -- +(1,1) node[anchor=west, rotate=45] {$0$}; 
\draw[thick] (90,90) ++(20,-20) +(-1,-1) -- +(1,1) node[anchor=west, rotate=45] {$40$}; 
\draw[thick] (90,90) ++(40,-40) +(-1,-1) -- +(1,1) node[anchor=west, rotate=45] {$80$}; 
\draw[thick] (90,90) ++(60,-60) +(-1,-1) -- +(1,1) node[anchor=west, rotate=45] {$120$}; 

\end{tikzpicture}
\caption{The `plus' markers on the $z$-axis indicate eigenvalues of a large random matrix $A'$.
A large Young tableau analogous to Figure \ref{fig:tableau}  (individual boxes were not shown for clarity). Shaded regions show boxes with numbers smaller than some thresholds.
The diagonal lines show the $z$-coordinates of the concave corners of the Young diagram
drawn with a thick line.}
\label{figure:law-of-large}
\end{figure}

Voiculescu initiated a highly non-commutative probability theory called \emph{free probability}.
One of its highlights is related to \emph{random matrices}.
We consider the following concrete problem.
Let 
$$ A = \begin{bmatrix} a_{11} & \cdots & a_{1n} \\
                       \vdots & \ddots & \vdots \\
                       a_{n1} & \cdots & a_{nn} 
       \end{bmatrix}$$ 
be an $n\times n$ random matrix, selected uniformly from the manifold of all hermitian matrices
with prescribed eigenvalues $x_1,\dots,x_n$. 
Let $1\leq m< n$; \emph{what can we say about the eigenvalues of the 
$m\times m$ upper-left corner}  
$$ A' = \begin{bmatrix} a_{11} & \cdots & a_{1m} \\
                       \vdots & \ddots & \vdots \\
                       a_{m1} & \cdots & a_{mm} 
       \end{bmatrix}?
$$
In the limit as $n\to\infty$ and $\frac{m}{n}$ converges to some limit, 
a kind of law of large numbers occurs and 
these eigenvalues with high probability concentrate around some limit distribution
depending only on the distribution of the eigenvalues $x_1,\dots,x_n$ of the big matrix.
As an illustration, we present the
results of a computer experiment for a large matrix $A$ 
with eigenvalues $-80,0,120$ (each with some high multiplicity which will be discussed later).
The eigenvalues of
the corner matrix $A'$ are shown as `plus' markers on the $z$-axis on Figure \ref{figure:law-of-large}. 
Some of these eigenvalues are degenerate and coincide with the three eigenvalues of the original matrix
(thick markers). 
The remaining eigenvalues occupy 
two intervals: one around $-40$ and one around $80$,
the empirical density of these eigenvalues turns out to be quite close to the asymptotic value for $n\to\infty$.

In order to explain this law of large numbers phenomenon 
we will use  \emph{free cumulants}. 
The basic idea is that
even though the random matrix $A$ is a complicated, multidimensional object, its corner entry $a_{11}$ is just
a complex-valued random variable which can be investigated by the (logarithm of) the Fourier-Laplace transform.
The fact that our random matrix $A$ has a large symmetry 
implies that this corner entry $a_{11}$ contains essentially all information
about $A$ which is necessary for asymptotic problems. 
The free cumulants
$R_k=R_k(A)$ of the random matrix $A$ (with $k\in\{1,2,\dots\}$) are defined as suitably normalized 
coefficients in the expansion 
$$ \log \E e^{t a_{11}} = t R_1 + \frac{t^2}{2 n} R_2 + \frac{t^3}{3 n^2} R_3 + \cdots.$$
The normalization was chosen in such a way that the free cumulants converge to finite values as $n\to\infty$.

Free cumulants contain the same information as the eigenvalues of the matrix but they are much more 
convenient for asymptotic problems. 
For example, the solution to the our problem of the eigenvalues of the corner $A'$ is given by 
\begin{equation}
\label{eq:random-matrix}
 R_{k+1}(A') = \left(\frac{m}{n}\right)^{k} R_{k+1}(A) 
\end{equation}
which is a direct consequence of the above definition of free cumulants.

\subsection{Free cumulants of a Young diagram}


It was observed by Biane that to an irreducible representation $\rho^\lambda$ of a symmetric group
one can associate a certain large matrix $\Gamma^\lambda$
which contains all information about $\rho^\lambda$.
The eigenvalues of this matrix are nicely related to the shape of the Young diagram: 
the $z$-coordinates (defined as $x-y$) of the concave corners of the 
Young diagram (indicated by the diagonal gray lines on Figure \ref{fig:youngA}) are the eigenvalues of this matrix.
For example, the matrix $\Gamma^\lambda$ 
associated to the Young diagram from Figure \ref{fig:youngA} has eigenvalues: $-2$, $0$, $3$
(each with some high multiplicity).
In particular, the eigenvalues of $\Gamma^{\lambda}$ and of $\Gamma^{s\lambda}$ are related to
each other by a simple scaling by factor $s$; compare Figures \ref{fig:youngA} and \ref{fig:youngB}.

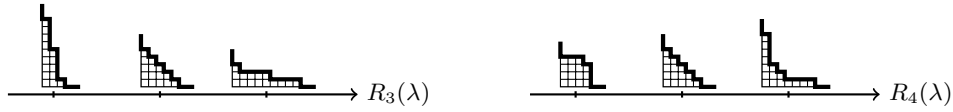
\begin{figure}[t]
\centering 
\begin{tikzpicture}[scale=0.2]
\small
\draw[->,thick] (-10,0) -- (13,0) node[anchor=west]{$R_3(\lambda)$};
\foreach \x in {-7,0,7}
                \draw[thick] (\x,6pt) -- (\x,-6pt);

\begin{scope}[scale=0.5,shift={(-2.5,1)}]
\draw (5,0) -- (0,0) -- (0,5); 
\draw[ultra thick] (7,0) -- (5,0) -- (5,1) -- (4,1) -- (4,2) -- (3,2) -- (3,3) -- (2,3) -- (2,4) -- (1,4) -- (1,5) -- (0,5) -- (0,7);
\clip (0,0) -- (5,0) -- (5,1) -- (4,1) -- (4,2) -- (3,2) -- (3,3) -- (2,3) -- (2,4) -- (1,4) -- (1,5) -- (0,5);
\draw (0,0) grid (5,5);
\end{scope}

\begin{scope}[scale=0.5,shift={(9.5,1)}]
\draw (9,0) -- (0,0) -- (0,3);
\draw[ultra thick] (11,0) -- (9,0) -- (9,1) -- (5,1) -- (5,2) -- (1,2) -- (1,3) -- (0,3) -- (0,5);
\clip (0,0) -- (9,0) -- (9,1) -- (5,1) -- (5,2) -- (1,2) -- (1,3) -- (0,3);
\draw (0,0) grid (9,5);
\end{scope}

\begin{scope}[scale=0.5,shift={(-15.5,1)}]
\draw(0,9) -- (0,0) -- (3,0);
\draw[ultra thick] (0,11) -- (0,9) -- (1,9) -- (1,5) -- (2,5) -- (2,1) -- (3,1) -- (3,0) -- (5,0);
\clip (0,0) -- (0,9) -- (1,9) -- (1,5) -- (2,5) -- (2,1) -- (3,1) -- (3,0);
\draw (0,0) grid (5,9);
\end{scope}
\end{tikzpicture}
\hfill
\begin{tikzpicture}[scale=0.2]
\small
\draw[->,thick] (-10,0) -- (13,0) node[anchor=west]{$R_4(\lambda)$};
\foreach \x in {-7,0,7}
                \draw[thick] (\x,6pt) -- (\x,-6pt);
\begin{scope}[scale=0.5,shift={(-2.5,1)}]
\draw (5,0) -- (0,0) -- (0,5);
\draw[ultra thick] (7,0) -- (5,0) -- (5,1) -- (4,1) -- (4,2) -- (3,2) -- (3,3) -- (2,3) -- (2,4) -- (1,4) -- (1,5) -- (0,5) -- (0,7);
\clip (0,0) -- (5,0) -- (5,1) -- (4,1) -- (4,2) -- (3,2) -- (3,3) -- (2,3) -- (2,4) -- (1,4) -- (1,5) -- (0,5);
\draw (0,0) grid (5,5);
\end{scope}

\begin{scope}[scale=0.5,shift={(10.5,1)}]
\draw (7,0) -- (0,0) -- (0,7);
\draw[ultra thick] (9,0) -- (7,0) -- (7,1) -- (3,1) -- (3,2) -- (2,2) -- (2,3) -- (1,3) -- (1,7) -- (0,7) -- (0,9); 
\clip (0,0) -- (7,0) -- (7,1) -- (3,1) -- (3,2) -- (2,2) -- (2,3) -- (1,3) -- (1,7) -- (0,7);
\draw (0,0) grid (7,7);
\end{scope}

\begin{scope}[scale=0.5,shift={(-16,1)}]
\draw (4,0) -- (0,0) -- (0,4);
\draw[ultra thick] (6,0) -- (4,0) -- (4,3) -- (3,3) -- (3,4) -- (0,4) -- (0,6); 
\clip (0,0) -- (4,0) -- (4,3) -- (3,3) -- (3,4) -- (0,4) -- (0,6);
\draw (0,0) grid (5,9);
\end{scope}

\end{tikzpicture}

\caption{
Intuitive meaning of the free cumulants $R_3$ and $R_4$ as parameters describing shape of a Young diagram.}
\label{fig:r34}
\end{figure}

Biane defined the free cumulants of the Young diagram $\lambda$
$$ R_k(\lambda):= R_k(\Gamma^\lambda) $$
as free cumulants of the corresponding matrix $\Gamma^\lambda$. 
The free cumulant $R_k$ is a homogeneous function of degree $k$ on the set of Young diagrams:
$$ R_k(s\lambda) = s^k R_k(\lambda),$$
in other words free cumulants are examples of the `nice' functions which we were looking for in Section \ref{subsec:homogeneous}.
From our perspective we can forget that free cumulants of Young diagrams have this long and interesting
history related to random matrices and simply treat them as convenient parameters describing the shape of the Young diagram and which can be calculated efficiently.
For example, for a Young diagram $\lambda$ with $n$ boxes (Figure \ref{fig:r34})
$$ R_3(\lambda)= 2 \iint\limits_{(x,y)\in\lambda} (x-y) \ dx\ dy; \qquad 
 R_4(\lambda)= 3 \iint\limits_{(x,y)\in\lambda} (x-y)^2 \ dx\ dy -\frac{3}{2} n.$$

\section{Kerov polynomials}
\label{sec:kerov}

It was proved by Kerov that free cumulants $R_k=R_k(\lambda)$ can be used for calculation of the
normalized characters $\Ch_k=\Ch_k(\lambda)$. For example,
\begin{align*}
\Ch_2 &= R_3, \\
\Ch_3 &= R_4 + R_2,   \\
\Ch_4 &= R_5 + 5R_3,    \\     
\Ch_5 &= R_6 + 15R_4 + 5R_2^2 + 8R_2.
\end{align*}
The right-hand sides are called \emph{Kerov polynomials}.
The reader can recognize that our guiding formula \eqref{eq:main} is among them.
We will discuss some interesting features of these polynomials.

\subsection{$\Ch_5\approx R_6$}
We evaluate the equality \eqref{eq:main} on the dilated diagram $s\lambda$.
Homogeneity of free cumulants implies that
$$ \Ch_5(s\lambda) = \underbrace{s^6 R_6(\lambda)}_{\text{degree $6$}} + \underbrace{15 s^4 R_4(\lambda)+ 5s^4 R_2^2(\lambda)}_{\text{degree $4$}}  + \underbrace{8 s^2R_2(\lambda)}_{\text{degree $2$}}. $$
In the limit $s\to\infty$ only the top-degree part really matters, therefore we can informally write
$ \Ch_5 \approx R_6 $.
It was shown by Biane that it is a general phenomenon: the value of the (normalized) irreducible character on the cycle $[k]$ is (asymptotically, for large Young diagrams) given by the free cumulant $R_{k+1}$ of the Young diagram:
$$ \Ch_k \approx R_{k+1}. $$
The left-hand side is the quantity which we wanted to understand because it is so fundamental for the representation
theory; the right-hand side can be efficiently calculated from the shape of the Young diagram.
This is a beautiful result and we will present one of its applications below.

\subsection{Biane's law of large numbers}
\label{subsec:biane}

Let us randomly select a Young tableau filling prescribed Young diagram $\lambda$ with $n$ boxes and let us remove
boxes with numbers bigger than some prescribed threshold $1\leq m<n$ (Figure \ref{fig:tableau}). 
\emph{What can we say about the shape of the resulting smaller Young diagram $\mu$?}
The results of computer experiments (Figure \ref{figure:law-of-large})
suggest that with high probability these new Young diagrams asymptotically concentrate around some smooth limiting shapes.

This problem can be reformulated in the language of the representation theory:
$\rho^\lambda$ is a representation of $\Sym{n}$; 
we consider its restriction $\rho^\lambda\downarrow_{\Sym{m}}$ to the subgroup $\Sym{m}$. This restriction
is usually a reducible representation; we decompose it into irreducible components and randomly select one of them,
say $\rho^\mu$. The distribution of the resulting random Young diagram $\mu$ has the same distribution as in the original problem. 

For large Young diagrams $\lambda$ and $\mu$ we can write
\begin{align*}R_{k+1}(\lambda) & \approx \Ch_k(\lambda) \approx n^k  \frac{\Tr \chi^\lambda([k])}{\Tr \chi^\lambda(e)},\\
\E R_{k+1}(\mu) & \approx \E \Ch_k(\mu) \approx m^k \E \frac{\Tr \chi^\mu([k])}{\Tr \chi^\mu(e)} 
=m^k \frac{\Tr \chi^\lambda([k])}{\Tr \chi^\lambda(e)}. \end{align*}
By comparing these two approximate equalities we conclude that for a \emph{typical} random
Young diagram $\mu$ we can expect that
\begin{equation}
\label{eq:compression-young}
 R_{k+1}(\mu) \approx \left(\frac{m}{n}\right)^k R_{k+1}(\lambda).  
\end{equation}
What a surprise! This equality has the same form as \eqref{eq:random-matrix} for the random matrices.
This is illustrated on Figure \ref{figure:law-of-large}: the diagonal lines indicate the $z$-coordinates
of the concave corners of the Young diagram $\mu$ drawn with a thick line while the `plus' markers indicate
eigenvalues of a corner of a large random matrix $A$ with the same eigenvalues 
as $\Gamma^\lambda$.
The parallelism between \eqref{eq:random-matrix} and \eqref{eq:compression-young} implies that 
(asymptotically) the density of the eigenvalues of the corner matrix $A'$ should match the 
density of eigenvalues of $\Gamma^\mu$ and thus the $z$-coordinates of 
the concave corners of the diagram $\mu$.
As one can see on Figure \ref{figure:law-of-large}, there is indeed a good match.

\subsection{Kerov positivity conjecture}
An interesting feature of the examples of Kerov polynomials presented at the beginning of 
Section \ref{sec:kerov} is that
their coefficients are non-negative integers. 
\emph{Kerov positivity conjecture} states that it is a general phenomenon.
The fact that these coefficients are integers followed easily from 
Kerov's construction but
their positivity was rather mysterious. 

Combinatorialists tend to believe that
\emph{if some reasonable integer numbers turn out to be non-negative, 
there should be a natural explanation by showing
that they are cardinalities of some interesting objects}. 
Following this line of thinking, positivity conjectures indicate that the object we 
are investigating might have some hidden underlying structure and 
thus such conjectures are very inspirational for the research.
This was also the case with Kerov conjecture; it initiated investigation of the 
characters of the symmetric groups with a new perspective. We will review it in the following.

\subsection{Maps}

\begin{figure}
\centering
\begin{tikzpicture}[scale=0.5,
white/.style={circle,draw=black,inner sep=4pt},
black/.style={circle,draw=black,fill=black,inner sep=4pt},
connection/.style={draw=black!80,black!80,auto}
]
\footnotesize

\begin{scope}
\clip (0,0) rectangle (10,10);
\draw (3,5) node (b1) [black,label=80:$R_4$] {};
\draw (b1) +(10,0) node (b1prim) [black] {};

\draw (b1) +(1,-3) node (b1-se) [white] {};
\draw (b1) +(-1,-3) node (b1-sw) [white] {};

\draw (8,8) node (b2) [black,label=-60:$R_2$] {};
\draw (b2) +(0,-10) node (b2prim) [black] {};
\draw (b2) +(-10,0) node (b2prim2) [black] {};

\draw (6,5) node (w1) [white] {};
\draw (w1) +(0,10) node (w1prim) {};

\draw (12,7) node (w2) [white] {};
\draw (w2) +(-10,0) node (w2prim) [white] {};

\draw[connection,pos=0.2] (b2) to node {$4$} node [swap] {} (w1prim);
\draw[connection,pos=0.666] (b2prim) to node {$4$} node [swap] {} (w1);

\draw[connection,pos=0.2] (b2) to node {$6$} node [swap] {} (w2);
\draw[connection,pos=0.7] (b2prim2) to node {$6$} node [swap] {} (w2prim);

\draw[connection] (b1) to node {$2$} node [swap] {} (b1-sw);

\draw[connection] (b1) to node {$3$} node [swap] {} (b1-se);

\draw[connection] (b1) to node {$5$} node [swap] {} (w1);

\draw[connection] (w1) to node {} node [swap] {$7$} (b2);

\draw[connection] (w2) to node {} node [swap] {$1$} (b1prim);
\draw[connection] (w2prim) to node {} node [swap] {$1$} (b1);

\end{scope}

\draw[very thick,decoration={
    markings,
    mark=at position 0.5 with {\arrow{>}}},
    postaction={decorate}]  
(0,0) -- (10,0);

\draw[very thick,decoration={
    markings,
    mark=at position 0.5 with {\arrow{>}}},
    postaction={decorate}]  
(0,10) -- (10,10)  ;

\draw[very thick,decoration={
    markings,
    mark=at position 0.5 with {\arrow{>>}}},
    postaction={decorate}]  
(0,0) -- (0,10);

\draw[very thick,decoration={
    markings,
    mark=at position 0.5 with {\arrow{>>}}},
    postaction={decorate}]  
(10,0) -- (10,10)  ;

\end{tikzpicture}
\caption{Map on the torus. The left side of the square should be glued to the right side,
as well as bottom to top, as indicated by arrows.}
\label{fig:map-kerov}
\end{figure}
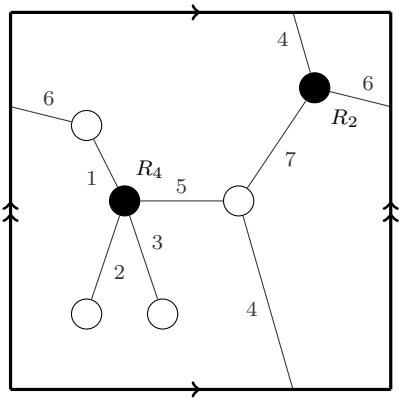

A \emph{map} is a graph drawn on an oriented surface (Figure \ref{fig:map-kerov}). 
We will consider only maps which are \emph{bipartite} 
(each vertex is either white or black, there are no edges between vertices of the same color),
\emph{unicellular} (if we remove the graph from the surface, the remaining part --- called \emph{cell} --- is homeomorphic to one disc), \emph{labeled} (the edges are labeled; if we go clockwise along the boundary of the cell and read every second label, whey will form the sequence $1,2,\dots,k$).

A map carries more information than just a graph, for example for each vertex it makes sense to speak about
the cyclic order of the incident edges.
We can encode this by a cycle from the permutation group $\Sym{k}$.
By merging such disjoint cycles corresponding to white vertices (respectively, black vertices) we obtain a permutation $\sigma_1$ (respectively, $\sigma_2$). For example, map presented on Figure \ref{fig:map-kerov} corresponds to 
$\sigma_1=(1,6)(2)(3)(4,7,5)$ and $\sigma_2=(1,2,3,5)(4,7,6)$. Permutations $\sigma_1$ and $\sigma_2$ contain the same
information as the original map. In particular, the structure of the cells of our map can be recovered from the 
product $\sigma_1 \sigma_2$; in our example
$$ \sigma_1 \sigma_2 = (1,2,3,\dots,7) = [7] \in \Sym{7}$$
has exactly one cycle which reflects the fact that our map is unicellular. 
Studying the maps is therefore equivalent to studying solutions of the equation
$$ \sigma_1 \sigma_2 = [k]\qquad \text{with }\sigma_1,\sigma_2\in \Sym{k},$$ 
but maps have an advantage related to their geometric and graph-theoretic flavor.

\subsection{Stanley character formula}

\begin{figure}
\centering
\subfloat[][]{\begin{tikzpicture}[scale=0.5,
white/.style={circle,draw=black,inner sep=4pt},
black/.style={circle,draw=black,fill=black,inner sep=4pt},
connection/.style={draw=black!80,black!80,auto}
]
\footnotesize

\begin{scope}
\clip (0,0) rectangle (10,10);

\draw (3.333,2.333) node (b1)    [black,label=90:$\Pi$] {};
\draw (b1) +(10,0) node (b1prim) [black] {};

\draw (7.666,6.666) node (b2)     [black,label=0:$\Sigma$] {};
\draw (b2) +(0,-10) node (b2prim) [black] {};

\draw (b2) +(-3,1) node (w2) [white,label=180:$W$] {};

\draw (6.666,3.333) node (w1) [white,label=45:$V$] {};
\draw (w1) +(-10,0) node (w1left) [white] {};
\draw (w1) +(0,10)  node (w1top)  [white] {};

\draw[connection] (b1) to node {$4$} node [swap] {} (w1);

\draw[connection] (b2) to node {$3$} node [swap] {} (w2);

\draw[connection,pos=0.333] (b2) to node {$2$} node [swap] {} (w1top);
\draw[connection,pos=0.666] (b2prim) to node {$2$} node [swap] {} (w1);

\draw[connection,pos=0.666] (b1prim) to node {$1$} node [swap] {} (w1);
\draw[connection,pos=0.333] (b1) to node {$1$} node [swap] {} (w1left);

\draw[connection,pos=0.666] (w1) to node {} node [swap] {$5$} (b2);

%
%
%
%
%
%
%
%

\end{scope}

\draw[very thick,decoration={
    markings,
    mark=at position 0.5 with {\arrow{>}}},
    postaction={decorate}]  
(0,0) -- (10,0);

\draw[very thick,decoration={
    markings,
    mark=at position 0.5 with {\arrow{>}}},
    postaction={decorate}]  
(0,10) -- (10,10)  ;

\draw[very thick,decoration={
    markings,
    mark=at position 0.5 with {\arrow{>>}}},
    postaction={decorate}]  
(0,0) -- (0,10);

\draw[very thick,decoration={
    markings,
    mark=at position 0.5 with {\arrow{>>}}},
    postaction={decorate}]  
(10,0) -- (10,10)  ;

\end{tikzpicture}
\label{subfig:map}}
\hfill
\subfloat[][]{
\begin{tikzpicture}[scale=1.2]
\begin{scope}

\draw[line width=5pt,black!20] (-0.2,0.5) -- (3.2,0.5);
\draw (3.2,0.5) node[anchor=west] {$\Sigma$};

\draw[line width=5pt,black!20] (-0.2,1.5) -- (1.2,1.5);
\draw (1.2,1.5) node[anchor=west] {$\Pi$};

\draw[line width=5pt,black!20] (2.5,-0.2) -- (2.5,1.2);
\draw (2.5,1.2) node[anchor=south] {$W$};

\draw[line width=5pt,black!20] (0.5,-0.2) -- (0.5,2.2);
\draw (0.5,2.2) node[anchor=south] {$V$};

\draw[ultra thick] (4.5,0) -- (3,0) -- (3,1) -- (1,1) -- (1,2) -- (0,2) -- (0,2.5); 
\clip (0,0) -- (3,0) -- (3,1) -- (1,1) -- (1,2) -- (0,2); 
\draw (0,0) grid (3,3);
\end{scope}
\draw (0.5,-0.2) node[anchor=north,text height=8pt] {$a$};
\draw (1.5,-0.2) node[anchor=north,text height=8pt] {$b$};
\draw (2.5,-0.2) node[anchor=north,text height=8pt] {$c$};

\draw (-0.2,0.5) node[anchor=east]  {$\alpha$};
\draw (-0.2,1.5) node[anchor=east] {$\beta$};

\draw(2.5,0.5) node {$3$};
\draw(0.5,0.5) node {$2,5$};
\draw(0.5,1.5) node {$1,4$};

\end{tikzpicture}
\label{subfig:embed}}

\caption{\protect\subref{subfig:map} map on the torus and 
\protect\subref{subfig:embed} an example of its embedding 
$F(\Sigma)=\alpha$, $F(\Pi)=\beta$, $F(V)=a$, $F(W)=c$.
$F(1)=F(4)=(a \beta)$, $F(2)=F(5)=(a \alpha)$, $F(3)=(c\alpha)$.
The columns of the Young diagram were indexed by Latin letters, the rows by Greek letters.}
\label{fig:embedding}
\end{figure}
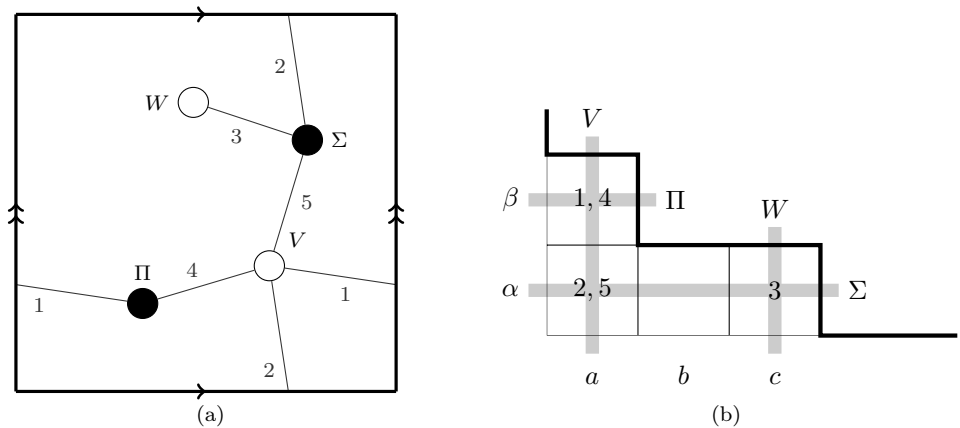

Attempts to prove Kerov conjecture have led in a natural way to discovery of
Stanley's formula for normalized characters:
$$ \Ch_k(\lambda)= \sum_{M} (-1)^{k-\#\text{white vertices}} N_M(\lambda), $$
where the sum runs over all maps $M$ with $k$ edges.
Above $N_M(\lambda)$ denotes the number of \emph{embeddings of the map $M$ to the
Young diagram $\lambda$} (Figure \ref{fig:embedding}).
An embedding is a function which maps white vertices to columns of $\lambda$, black vertices to rows
of $\lambda$, edges to boxes of $\lambda$. We also require that an embedding preserves the incidence,
i.e.~a vertex $V$ and an incident edge $E$ should be mapped to a row or column $F(V)$ which contains the box
$F(E)$. 

Stanley's formula is a perfect tool for studying asymptotics of characters of symmetric groups in various scalings.
It was also the tool which was essential in the proof of Kerov positivity conjecture (which will be discussed below).

\subsection{Genus expansion}
The function $\lambda\mapsto N_M(\lambda)$ is homogeneous which explains why Stanley's formula is
a perfect tool for our purposes. The degree of this homogeneous function
$$ \operatorname{deg} N_M = k+1 - 2 \operatorname{genus}(M) $$
is directly related to the \emph{genus} of the surface on which map $M$ is drawn. Thus the \emph{planar maps}
(which can be drawn on the sphere, genus equal to zero) have maximal possible degree and asymptotically
have the biggest contribution.
This kind of \emph{genus expansion} where to combinatorial summands one can associate a surface
which determines the asymptotics is very common in the asymptotic representation theory as well as in the random matrix theory. Genus expansion and, in particular, the special role of planar maps explains why combinatorics of 
free cumulants (originally formulated by Speicher in terms of \emph{non-crossing partitions} which are set partitions that can be drawn on a sphere) is so useful.

\subsection{Proof of Kerov's conjecture and combinatorial interpretation of Kerov polynomials}
\label{subsec:Kerov-interpretation}

The coefficient standing at monomial $R_{i_1} \cdots R_{i_l}$ in Kerov polynomial
$\Ch_k$ turns out to be the number of maps with $k$ edges with black vertices 
decorated by $R_{i_1},\dots,R_{i_l}$ (Figure \ref{fig:map-kerov}) such that the following
\emph{transportation problem} has a solution. We imagine that each white vertex is a factory producing a unit of some liquid, each black vertex decorated by $R_i$ is a consumer demanding $i-1$ units of this liquid and the edges of the map are one-way pipes which can transport the liquid only from white to black vertices. We require that it is
possible to arrange the amount of the liquid in each pipe in such a way that \emph{each pipe transports a strictly positive amount}.
The map on Figure \ref{fig:map-kerov} fulfills this condition. 

The fact that we require a \emph{strictly} positive solution is quite unusual for such transportation problems and has some interesting consequences. Firstly, the classical criterion 
(given by \emph{Hall's marriage theorem}) for checking whether such a transportation problem has a solution has to 
be changed. 

Secondly, this \emph{strict positivity} requirement restricts the maps which
could contribute to the coefficients of Kerov polynomials, namely such a map cannot contain a disconnecting edge
except for edges leading to white leaves (Figure \ref{fig:map-kerov}).
This is quite a strong restriction; for example 
the number of such maps with a fixed genus grows only polynomially with the size of the map while the number of all maps with fixed genus grows exponentially.
This implies that the coefficients of Kerov polynomials for $\Ch_k$ corresponding to fixed genus grow relatively slowly (polynomially) with $k$.
It should be compared with analogues of Kerov polynomials in which instead of free cumulants we use some other quantities describing the shape of the Young diagram; in the latter case the growth of the integer 
coefficients with $k$ is usually
exponential. This is an indication that Kerov polynomials contain relatively small amount of information, they have small complexity and thus free cumulants are \emph{the right} quantities for
studying asymptotics of characters.

\subsection{Gaussian fluctuations}
The definition of the normalized characters
can be easily adapted to more complicated conjugacy classes; for example we denote by $\Ch_{k,l}(\lambda)$ the normalized character on a pair of cycles
of lengths $k$ and $l$.
Following our guide \eqref{eq:main} we can write an analogue of Kerov polynomials,
for example
$$ \Ch_{3,2} = R_3 R_4 - 5 R_2 R_3 - 6 R_5 - 18 R_3. $$
As one can see, the analogue of Kerov's positivity conjecture does not hold true.
This is an indication that \emph{characters $\Ch_{k,l}$ are not the right quantities}.

It turns out that it is much better to study a kind of \emph{covariance}
$$ \operatorname{Cov}(\Ch_k,\Ch_l):= \Ch_{k,l} - \Ch_k \Ch_l $$
which measures how the character on two disjoint cycles differs from the product of the characters on each cycle separately. One can consider the corresponding Kerov polynomials, for example
$$ \operatorname{Cov}(\Ch_3,\Ch_2):= \Ch_{3,2} - \Ch_3 \Ch_2 =  -\big( 6 R_2 R_3 + 6 R_5 + 18 R_3 \big); $$
apart from the global change of the sign all coefficients are again non-negative integers which is an indication that
\emph{such a covariance is the right  quantity}. Indeed, the combinatorial interpretation of the coefficients
of Kerov polynomials from Section \ref{subsec:Kerov-interpretation} holds true after some simple adjustments,
including the requirement that we consider only \emph{connected} maps (for unicellular maps this was automatic).

A degree of a function $F$ on Young diagrams is defined as the degree of the polynomial $s\mapsto F(s\lambda)$.
The \emph{connectivity} requirement in the combinatorial interpretation of Kerov polynomials
influences the topology of the maps which we count; hence the degree of the covariance
$\operatorname{Cov}(\Ch_k,\Ch_l)$
is smaller than the degrees of individual summands 
$\Ch_k \Ch_l$ and $\Ch_{k,l}$ which will have interesting consequences.

If the number of cycles is bigger, instead of covariance one should consider a \emph{cumulant} $k(\Ch_k,\Ch_l,\dots,\Ch_s)$
which measures in a more refined way how much the character $\Ch_{k,l,\dots,s}$ differs from products of characters
with simpler cycle structure. All the above mentioned results hold true also in this more general setup. 
The fact that the degree of the cumulant $k(\Ch_k,\Ch_l,\dots,\Ch_s)$ is much smaller than the sum of degrees of 
the individual factors implies that (if proper normalization is chosen)   
$\Ch_2,\Ch_3,\dots$ regarded in a rather abstract way as random variables are asymptotically Gaussian. 

More specifically, this implies that a generalization of 
\emph{Kerov's Central Limit Theorem} holds true: for a wide class of reducible representations of the symmetric groups
if we randomly select an irreducible component $\rho^\mu$ (like in Section \ref{subsec:biane}), the fluctuations of the shape of $\mu$
will be asymptotically Gaussian. This is yet another application of Kerov polynomials and their combinatorics.

\section{Open problems}

As a rule, open problems are much more interesting than the solved ones. 
Fortunately, there are still several mysteries concerning Kerov polynomials and related
objects. 
We will review some of these open problems.
Just like Kerov's conjecture, they are also related to positivity, thus they
hint at some unexpected hidden combinatorial structures and we hope that investigation of them will
be as profitable as the investigation of Kerov conjecture was.

\subsection{Goulden-Rattan character polynomials}

Goulden-Rattan polynomials express \emph{the difference} $\Ch_{k}-R_{k+1}$ as a polynomial in $C_2,C_3,\dots$
(where  $C_2(\lambda),C_3(\lambda),\dots$  are some quantities, related to free cumulants, describing the shape of a Young diagram).
For example
$$ \Ch_7 - R_8 = 14 C_6 + \frac{469}{3} C_4 + \frac{203}{3}
C_2^2 + 180 C_2. $$
These polynomials are less complicated then the corresponding Kerov polynomials 
(i.e., contain a smaller number of summands), which suggests that $(C_k)$ are better, \emph{more fundamental},
for understanding the deviation from the approximation $\Ch_{k}\approx R_{k+1}$. Furthermore, the coefficients
of these polynomials seem to be positive rational numbers with relatively small denominators. 
It would be more difficult to find a combinatorial interpretation of positive numbers which are not integers,
nevertheless this possibility is very tempting.

\subsection{Kerov polynomials for Jack characters}
Lassalle observed that just like normalized characters describe the dual combinatorics
of representations of symmetric groups, it is possible to consider the dual combinatorics of \emph{Jack polynomials}.
In this way one can obtain a quite natural deformation of the characters of the symmetric group
with an additional parameter $\gamma$. Also in this more general case it is possible to find Kerov polynomials, for example
$$    \Ch^{(\gamma)}_{4} = R_5 + 6\gamma R_4 + \gamma R_2^2 + (5 + 11\gamma^2)R_3 + (7\gamma +
6\gamma^3)R_2.$$
As the reader can see, the coefficients of these more general Kerov polynomials also seem to be non-negative integers.

\section{Further reading}

Due to lack of space we will refer mostly to overview articles.
Ref.~\cite{Novak'Sniady2011} is a haiku-style introduction to free cumulants.
Ref.~\cite{Biane2002} is an overview of combinatorics of free cumulants and their applications
to random matrices and representation theory.
A lengthy introduction to \cite{DolkegaF'eray'Sniady2010} gives an overview of  Kerov polynomials.
Ref.~\cite{F'eray'Sniady2011} gives more details on Stanley's character formula and its applications
to asymptotics of characters. Ref.~\cite{'Sniady2006} gives more details on Gaussian fluctuations of Young diagrams.


\end{document}